\documentclass{article}

\usepackage[latin1]{inputenc}
\usepackage{subfigure}

\usepackage{amsfonts}
\usepackage{amssymb}
\usepackage{amsmath}
\usepackage{url}

\newtheorem{lemma}{Lemma}
\newtheorem{example}{Example}
\newtheorem{theorem}{Theorem}
\newtheorem{proposition}{Proposition}
\newtheorem{remark}{Remark}
\newtheorem{conjecture}{Conjecture}
\newtheorem{corollary}{Corollary}

\def\fff{{\mathbb{F}}}

\def\qqq{\mathbb{Q}}

\def\ccc{\mathbb{C}}

\def\pf{{\bf Proof}:\ }
\def\qed{$\Box$}
\def\wt{{\rm{wt}}}

\def\wt{{\rm{wt}}}
\def\sage{{\tt{SAGE}}}

\author{David Joyner}

\title{On quadratic residue codes and hyperelliptic curves}

\begin{document}
\maketitle
\begin{abstract}

For an odd prime $p$ and each non-empty subset $S\subset GF(p)$, consider 
the hyperelliptic curve $X_S$ defined by $y^2=f_S(x)$, where 
$f_S(x) = \prod_{a\in S}(x-a)$.  
Using a connection between binary quadratic residue codes and hyperelliptic 
curves over $GF(p)$, this paper investigates how coding theory bounds give rise to
bounds such as the following example:
for all sufficiently large primes $p$ there exists a 
subset $S\subset GF(p)$ for which the bound $|X_S(GF(p))| > 1.39p$
holds. We also use the quasi-quadratic residue codes defined below to
construct an example of a formally self-dual optimal code whose
zeta function does not satisfy the ``Riemann hypothesis.''

\end{abstract}

A long standing problem has been to develop ``good''
binary linear codes to be used for error-correction.  
This paper investigates in some detail an attack on this problem
using a connection between quadratic residue codes and hyperelliptic 
curves. Codes with this kind of relationship 
have been investigated in Helleseth \cite{H},
Bazzi-Mitter \cite{BM}, Voloch \cite{V1}, and Helleseth-Voloch \cite{HV}. 
This rest of this introduction is devoted to explaining 
in more detail the ideas discussed in later sections.

Let $\fff=GF(2)$ be the field with two elements
and $C\subset \fff^n$ denote a binary block code
of length $n$. 
For any two $\mathbf{x},\mathbf{y}\in{\fff}^{n}$,
let $d(\mathbf{x},\mathbf{y})$ denote the {\it Hamming metric}:

\begin{equation}
d(\mathbf{x},\mathbf{y})=|\{1\leq i\leq n\,\,|\,\,
x_{i}\not=y_{i}\}|.\end{equation}
The {\it weight} $\wt(\mathbf{x})$ of $\mathbf{x}$ is the number
of non-zero entries of $\mathbf{x}$.  
The smallest weight of any non-zero codeword 
is denoted $d$ -  the minimum distance if $C$ is linear. 
When $C$ is linear, denote the dimension of $C$ by $k$ and
call $C$ an $[n,k,d]_2$-code.

Denoting the volume of a Hamming sphere of radius $r$ in $\fff^n$
by $V(n,r)$,
%
%\[
%V(n,r)=\sum_{i=0}^r \left(\begin{array}{c} n\\ i\end{array}\right).
%\] 
the binary version of the {\it Gilbert-Varshamov bound} asserts that
(given $n$ and $d$) there is an $[n,k,d]_2$ code $C$ satisfying
$k\geq \log_2(\frac{2^n}{V(n,d-1)})$ \cite{HP}.  

\begin{conjecture}
(Goppa's conjecture \cite{JV},\cite{G}) 
The binary version of the Gilbert-Varshamov bound is asymptotically exact.
\end{conjecture}

For each odd prime $p>5$, a QQR code\footnote{This code
is defined in \S \ref{sec:BM} below.} 
is a linear code of length $2p$.
Like the quadratic residue codes, the length and dimension are
easy to determine but the minimum distance is more
mysterious. In fact, the weight of each codeword can be explicitly
computed in terms of the number of solutions in integers mod $p$
to a certain type of (``hyperelliptic'') polynomial 
equation. 
To explain the results better, some more notation is needed.

For our purposes, a {\it hyperelliptic curve} $X$ over $GF(p)$ is a 
polynomial equation of 
the form $y^2=h(x)$, where $h(x)$ is a polynomial with coefficients in 
$GF(p)$ with distinct roots\footnote{This overly simplistic definition
brings to mind the famous Felix Klein quote:
``Everyone knows what a curve is, until he has studied enough 
mathematics to become confused through the countless number 
of possible exceptions.'' Please see Tsafsman-Vladut \cite{TV} or
Schmidt \cite{Sc} for a rigorous treatment.}. 
The number of solutions to $y^2=h(x)$
mod $p$, plus the number of ``points at infinity'' on $X$, 
will be denoted $|X(GF(p))|$. This quantity can be related to a sum
of Legendre characters (see Proposition \ref{prop:hyperelliptic} below), 
thanks to classical 
work of Artin, Hasse, and Weil. This formula yields good estimates for
$|X(GF(p))|$ in many cases (especially when $p$ is large
compared to the degree of $h$). A long-standing problem has been to 
improve on the trivial estimate when $p$ is small compared to the 
degree of $h$. It turns out the work of Tarnanen \cite{T}
easily yields some non-trivial information on this problem
(see for example Lemma \ref{lemma:tarn} below), but the results given 
here improve upon this.

For each non-empty subset $S\subset GF(p)$, consider 
the hyperelliptic curve $X_S$ defined by $y^2=f_S(x)$, where 
$f_S(x) = \prod_{a\in S}(x-a)$.  
Let $B(c,p)$ be the statement: {\it For all subsets $S\subset GF(p)$, 
$|X_S(GF(p))| \leq  c\cdot p$ holds.} Note that $B(2,p)$ is
trivially true, so the statement $B(2-\epsilon,p)$, for some fixed 
$\epsilon >0$, might not be horribly unreasonable.

\begin{conjecture}
(``Bazzi-Mitter conjecture'' \cite{BM})
There is a $c\in (0,2)$ such that, for an infinite number of primes 
$p$ the statement $B(c,p)$ holds.
\end{conjecture}

It is remarkable that these two conjectures are related.
In fact, using QQR codes we show that if, for an infinite number of primes 
$p$ with $p\equiv 1\pmod 4$, $B(1.77,p)$ holds then 
Goppa's conjecture is false.
Although this is a new result, it turns out that it is
an easy consequence of the QQR construction given in \cite{BM} 
if you think about things in the right way. Using LQR 
codes\footnote{These codes will be defined in \S \ref{sec:LQR} below.}
we will remove the condition $p\equiv 1\pmod 4$ at a cost of slightly
weakening the constant $1.77$ (see Corollary \ref{cor:main}).

The spectrum and Duursma zeta function of these QQR codes is discussed in 
Section \ref{sec:duursma} below and 
some examples are given (with the help of the software
package {\tt SAGE} \cite{S}). 
We show that the analog of the Riemann hypothesis for the
zeta function of an optimal formally self-dual code is false 
using the family of codes
constructed in \S 2. The section ends with
some intriguing conjectures.

We close this introduction with a few open questions which, on the basis
of this result, seem natural.

{\it Question 1}: For each prime $p>5$ is there an 
effectively computable subset $S\subset GF(p)$ such that
$|X_S(GF(p))|$ is ``large''?

Here ``large'' is left vague but what is intended is
some quantity which is unusual. By Weil's estimate
(valid for ``small''-sized subsets $S$), we could expect
about $p$ points to belong to $|X_S(GF(p))|$. Thus ``large''
could mean, say, $>c\cdot p$, for some fixed $c>1$.

The next question is a strong version of the 
Bazzi-Mitter conjecture.

{\it Question 2}: Does there exist a $c<2$ such that, 
for all sufficiently large $p$ and all $S\subset GF(p)$, we have
$|X_S(GF(p))| < c\cdot p$?

In the direction of these questions, for
Question 1, a coding theory bound 
of McEliese-Rumsey-Rodemich-Welsh allows one to establish
the following result (see Theorem \ref{cor:main}):
{\it There exists a constant $p_0$ having the following
property: if $p\equiv 1 \pmod{4}$ and $p>p_0$ then there exists a 
subset $S\subset GF(p)$ for which the bound 
$|X_S(GF(p))| > 1.62p$ holds}\footnote{Moreover, we can remove the hypothesis 
$p\equiv 1 \pmod{4}$ if we assume Conjecture \ref{conj:CNQ}.}.
Unfortunately, the method of proof gives no clue how 
to compute $p_0$ or $S$.  
Using the theory of long quadratic-residue codes,
we prove the following lower bound (Theorem \ref{thrm:lqr}):
For all $p>p_0$ there exists a 
subset $S\subset GF(p)$ for which the bound 
$|X_S(GF(p))| > 1.39p$ holds. Again, we do not know what 
$p_0$ or $S$ is. 
%However, the theory of self-dual codes also provides the following lower
%bound for all $S$ (Corollary \ref{cor:lowerbound}):
%{\it If $p\equiv -1 \pmod 8$ then for all non-empty subsets $S\subset GF(p)$,
%the bound $|X_S(GF(p))| > \sqrt{p}$ holds}. 

Finally, Felipe Voloch \cite{V2} has kindly allowed the 
author to include some interesting 
explicit constructions (which do not use any theory of
error-correcting codes) in this paper (see \S \ref{sec:voloch} below).
First, he shows the following result:
{\it If $p\equiv 1\pmod 8$ then there exists an effectively computable
subset $S\subset GF(p)$ for which the bound $|X_S(GF(p))| > 1.5p$
holds.} A similar result holds for $p\equiv 3,7\pmod 8$.
Second, he gives a construction which answers Question 2 in the negative. 

\section{Cyclotomic arithmetic mod $2$}
\label{sec:1}

Let $R = \fff [x]/(x^p - 1)$ and $r_S \in R$ denotes the
polynomial 

\[
r_S(x) =\sum_{i\in S} x^i,
\]
where $S\subseteq GF(p)$. By convention, if $S=\emptyset$ is the 
empty set, $r_S=0$. We define the {\it weight} of $r_S$, denoted 
$\wt(r_S)$, to be the cardinality $|S|$. 
(In other words, identify in the obvious way
each $r_S$ with an element of $\fff^{p}$ and define the 
weight of $r_S$ to be the Hamming weight of the 
associated vector.).
For the set $Q$ of quadratic residues in $GF(p)^\times$ and the set 
$N$ of non-quadratic residues in $GF(p)^\times$, we have 
$\wt(r_Q)=\wt(r_N)=(p-1)/2$.
Note that $r_S^2=r_{2S}$, where $2S$ is the set of elements 
$2s\in GF(p)$, for $s\in S$. Using this fact and the 
quadratic reciprocity law, one can easily show that 
the following are equivalent:

\begin{itemize}
\item
$r_Q^2=r_Q$,
\item 
$2\in Q$ 
\item
$p\equiv \pm 1\pmod 8$.
\end{itemize}
Moreover, if $2\in N$ then $r_Q^2=r_N$.

Let $S,S_1,S_2, S_1'$ denote subsets of $GF(p)$, with 
$S_1\cap S_1'=\emptyset$, and let $S^c=GF(p)-S$ denote the complement.  
For $a\in GF(p)$, let

\[
H(S_1,S_2,a) =\{(s_1,s_2)\in S_1\times S_2\ |\ s_1+s_2 \equiv a\pmod p\}.
\]
In particular,
\begin{itemize}
\item
$H(S_1,S_2,a)=H(S_2,S_1,a)$,
\item
there is a natural bijection $H(GF(p),S,a)\cong S$,
\item
if $S_1\cap S_1'=\emptyset$ then 
$H(S_1,S_2,a)+H(S_1',S_2,a)=H(S_1+S_1',S_2,a)$.
\end{itemize}
Let

\[
h(S_1,S_2,a)=|H(S_1,S_2,a)| \pmod 2.
\]
Adding $|H(S_1,S_2,a)|+|H(S_1^c,S_2,a)|=|S_2|$
to $|H(S_1^c,S_2^c,a)|+|H(S_1^c,S_2,a)|=|S_1^c|$, we obtain 

\begin{equation}
\label{eqn:n-sum}
h(S_1,S_2,a)\equiv h(S_1^c,S_2^c,a)+|S_1^c|+|S_2|\pmod 2 .
\end{equation}
From the definition of $r_S$,

\[
r_{S_1}(x)r_{S_2}(x)=\sum_{a\in GF(p)} h(S_1,S_2,a)x^a
\]
in the ring $R$. Let $*:R\rightarrow R$ denote the
involution defined by $(r_S)^*=r_{S^c}=r_S+r_{GF(p)}$.
We shall see below that this is not an algebra involution.

\begin{lemma}
\label{lemma:prod}
For all $S_1,S_2\subset GF(p)$, we have
\begin{itemize}
\item
$|S_1|$ odd, $|S_2|$ even:
$r_{S_1}r_{S_2}=r_{S_1}^*r_{S_2}^*$ has even weight.

\item
$|S_1|$ even, $|S_2|$ even:
$(r_{S_1}r_{S_2})^*=r_{S_1}^*r_{S_2}^*$ has even weight.

\item
$|S_1|$ even, $|S_2|$ odd:
$r_{S_1}r_{S_2}=r_{S_1}^*r_{S_2}^*$ has even weight.

\item
$|S_1|$ odd, $|S_2|$ odd:
$(r_{S_1}r_{S_2})^*=r_{S_1}^*r_{S_2}^*$ has odd weight.

\end{itemize}
\end{lemma}

This lemma follows from the discussion above by a straightforward
argument.

Note that $R_{even}=\{r_S\ |\ |S| \ {\rm even}\}$,
is a subring of $R$ and, by the previous lemma, 
$*$ is an algebra involution on $R_{even}$.

%\begin{remark}
%The following seem to be true, based 
%on computer computations using \sage : 
%
%\[
%r_N(x)r_Q(x) =
%\left\{
%\begin{array}{ll}
%0,& p\equiv 1\pmod{8},\\
%1,& p\equiv 3\pmod{8},\\
%x+x^2+...+x^{p-1},& p\equiv 5\pmod{8},\\
%1+x+x^2+...+x^{p-1},& p\equiv 7\pmod{8}.
%\end{array}
%\right.
%\]
%These are not really needed, and I do have a proof of them, but
%they may be of interest for a reader wishing to 
%explore further the codes introduced in the next section.
%NOTE: Both sides are invariant under $x\to x^r$, so
%for $p\equiv 1\pmod{4}$, $r_N(x)r_Q(x)$ must be
%$0$ or $x+x^2+...+x^{p-1}$ for parity reasons.
%Likewise, for $p\equiv 3\pmod{4}$, $r_N(x)r_Q(x)$ must be
%$1$ or $1+x+x^2+...+x^{p-1}$. The quadratic reciprocity thrm and
%$h(N,Q,a)=h(N,Q,ra)$ might distinguish these cases.
%\end{remark}

\section{QQR Codes}
\label{sec:BM}

These are some observations on the interesting paper by Bazzi and Mitter 
\cite{BM}. We shall need to remove the 
assumption $p\equiv 3\pmod 8$ (which they make in their paper) below.

\vskip .2in

If $S \subseteq GF(p)$, let 
$f_S(x) = \prod_{a\in S}(x - a) \in GF(p)[x]$. 
Let $\chi = (\frac{\ }{p})$ be the quadratic
residue character, which is $1$ on the 
quadratic residues $Q\subset GF(p)^\times$, 
$-1$ on the quadratic non-residues $N\subset GF(p)^\times$, 
and is $0$ at $0\in GF(p)$.

Define 
\[
C_{NQ}=\{(r_N r_S, r_Q r_S) \ |\  S \subseteq GF(p)\}, 
\]
where $N,Q$ are as above. (We identify in the obvious way
each pair $(r_N r_S, r_Q r_S)$ with an element of $\fff^{2p}$.
In particular, when $S$ is the empty set, $(r_N r_S, r_Q r_S)$ 
is associated with the the zero vector in $\fff^{2p}$.) 
We call this a {\it QQR code} (or a {\it quasi-quadratic residue code}). 
These are binary linear codes of length $2p$ and dimension 

\[
k=
\left\{
\begin{array}{ll}
p, & {\rm if}\ p\equiv 3\pmod 4,\\
p-1, &{\rm if}\ p\equiv 1\pmod 4.
\end{array}
\right.
\]
This code has no codewords of odd weight, for parity reasons,
by Lemma \ref{lemma:prod}.

\begin{remark}
\label{remark:sqrt}
If $p\equiv \pm 1\pmod 8$ then $C_{NQ}$ ``contains'' 
a binary quadratic residue code. For such primes $p$, the 
minimum distance satisfies the well-known square-root
lower bound, $d\geq \sqrt{p}$.
\end{remark}

%Recall a linear code $C\subset \fff^{n}$ is {\it self-orthogonal} if
%$C$ is a subcode of its dual code 
%$C^\perp = \{v\in \fff^n\ |\ v\cdot c=0,\ \forall c \in C\}$.

Based on computations using \sage, the following statement is likely
to be true.

\begin{conjecture} 
\label{conj:CNQ}
For $p\equiv 1 \pmod{4}$, the associated QQR code and its dual satisfy:
$C_{NQ} \oplus C_{NQ}^\perp = \fff^{2p}$, where $\oplus$ stands for the 
direct product (so, in particular, $C_{NQ} \cap C_{NQ}^\perp = 
\{ \mathbf{0}\}$).
If $p\equiv 3 \pmod 4$ then the associated QQR code is self-dual:
$C_{NQ}^\perp = C_{NQ}$.
\end{conjecture}

The self-dual binary codes have useful upper bounds on their 
minimum distance (for example, the Sloane-Mallows bound
Theorem 9.3.5 in \cite{HP}). 
Combining this with the lower bound mentioned above, we have
the following result.

\begin{lemma}
\label{lemma:mindist}
Assume Conjecture \ref{conj:CNQ}.
If $p\equiv 3 \pmod 4$ then 
\[
d\leq 4\cdot [p/12]+6.
\]
If $p\equiv -1 \pmod 8$ then 
\[
\sqrt{p}\leq 
d\leq 4\cdot [p/12]+6.
\]
\end{lemma}

Note that these upper bounds (in the cases they are valid) 
are better than the asymptotic bounds of 
McEliese-Rumsey-Rodemich-Welsh for rate $1/2$ codes.

\begin{example}
The following computations were done with the help of \sage .
When $p=5$, $C_{NQ}$ has weight distribution
\[
[1, 0, 0, 0, 5, 0, 10, 0, 0, 0, 0].
\]
When $p=7$, $C_{NQ}$ has weight distribution
\[
 [1, 0, 0, 0, 14, 0, 49, 0, 49, 0, 14, 0, 0, 0, 1].
\]
When $p=11$, $C_{NQ}$ has weight distribution
\[
 [1, 0, 0, 0, 0, 0, 77, 0, 330, 0, 616, 0, 616, 0, 330, 0, 77, 0, 0, 0, 0, 0, 1].
\]
When $p=13$, $C_{NQ}$ has weight distribution
\[
[1, 0, 0, 0, 0, 0, 0, 0, 273,
 0, 598, 0, 1105, 0, 1300, 0,
 598, 0, 182, 0, 39, 0, 0, 0, 0, 0, 0].
\]

\end{example}

The following well-known result\footnote{See for example Weil \cite{W} or
Schmidt \cite{Sc}, Lemma 2.11.2.} shall be used to estimate the
weights of codewords of QQR codes.

\begin{proposition}
(Artin, Hasse, Weil)
\label{prop:hyperelliptic}
Assume $S\subset GF(p)$ is non-empty.
\begin{itemize}
\item
$|S|$ even: 
%The $\sum_{a\in GF(p)} \chi(f_S(a))$
%is equal to $-p-2$ plus the number of $GF(p)$-rational points on the
%(smooth projective model of the)
%hyperelliptic curve $X_S:\ y^2=f_{S}(x)$. In other words,

\[
\sum_{a\in GF(p)} \chi(f_S(a))
=-p-2+|X_S(GF(p))|.
\]

\item
$|S|$ odd: 
%The $\sum_{a\in GF(p)} \chi(f_S(a))$
%is equal to $-p-1$ plus the number of $GF(p)$-rational points on the
%(smooth projective model of the)
%hyperelliptic curve $X_S:\ y^2=f_{S}(x)$. In other words,

\[
\sum_{a\in GF(p)} \chi(f_S(a))
=-p-1+|X_S(GF(p))|.
\]

\item
$|S|$ odd: The genus of the (smooth projective model of the)
curve $y^2=f_S(x)$ is $g=\frac{|S|-1}{2}$ and

\[
\big{|}\sum_{a\in GF(p)} \chi(f_S(a))\big{|}
\leq (|S|-1)p^{1/2}+1.
\]

\item
$|S|$ even: The genus of the (smooth projective model of the)
curve $y^2=f_S(x)$ is $g=\frac{|S|-2}{2}$ and

\[
\big{|}\sum_{a\in GF(p)} \chi(f_S(a))\big{|}
\leq (|S|-2)p^{1/2}+1.
\]

\end{itemize}

\end{proposition}

Obviously, the last two estimates are only non-trivial
for $S$ ``small'' (e.g., $|S|<p^{1/2}$).

\begin{lemma}
\label{lemma:tarn}
(Tarnanen \cite{T}, Theorem 1)
Fix $\tau$, $0.39<\tau<1$. For all sufficiently large $p$,
the following statement is false: For all subsets $S\subset GF(p)$
with $|S|\leq \tau p$, we have $0.42p < |X_S(GF(p))|< 1.42p$.
\end{lemma}

\begin{remark}

\noindent
(1) Here the meaning of ``sufficiently large'' is hard to make
precise. The results of Tarnanen are actually asymptotic
(as $p\rightarrow \infty$), so we can simply say 
that the negation of part (1) of this Lemma contradicts 
Theorem 1 in \cite{T}.

\noindent
(2) This Lemma does not seem to imply 
``$B(1.42,p)$ is false, for sufficiently large $p$''
(so Theorem \ref{thrm:lqr} below is a new result), though it 
would if the condition $0.42p < |X_S(GF(p))|$ could be eliminated.
Also of interest is the statement about character sums in 
Theorem 1 of Stepanov \cite{St}.
\end{remark}

\pf
This is an immediate consequence of the Proposition above 
and Theorem 1 in \cite{T}. 
\qed

\begin{lemma}
(Bazzi-Mitter \cite{BM}, Proposition 3.3)
Assume $2$ and $-1$ are quadratic non-residues mod $p$ (i.e. 
$p \equiv 3 \pmod 8$).

If $\mathbf{c}=(r_N r_S, r_Q r_S)$ is a nonzero
codeword of the $[2p,p]$ binary code $C_{NQ}$ then the 
weight of this codeword can be expressed
in terms of a character sum as
\[
\wt(\mathbf{c})=p - \sum_{a\in GF(p)} \chi(f_S(a)),
\]
if $|S|$ is even, and
\[
\wt(\mathbf{c})=p +  \sum_{a\in GF(p)} \chi(f_{S^c}(a)),
\]
if $|S|$ is odd. 
\end{lemma}

In fact, looking carefully at their proof, one finds the following 
result.

\begin{proposition}
\label{prop:main}
Let $\mathbf{c}=(r_N r_S, r_Q r_S)$ be a nonzero
codeword of $C_{NQ}$. 
\begin{itemize}
\item[(a)] If $|S|$ is even
\[
\wt(\mathbf{c})=p - \sum_{a\in GF(p)} \chi(f_S(a))
=2p+2-|X_S(GF(p))|.
\]
 
\item[(b)] If $|S|$ is odd and $p\equiv 1\pmod 4$ then the weight is

\[
\wt(\mathbf{c})=p - \sum_{a\in GF(p)} \chi(f_{S^c}(a))
=2p+2-|X_{S^c}(GF(p))|.
\] 

\item[(c)] If $|S|$ is odd and $p\equiv 3\pmod 4$ then

\[
\wt(\mathbf{c})=p +  \sum_{a\in GF(p)} \chi(f_{S^c}(a))
=|X_{S^c}(GF(p))|-2.
\]
\end{itemize}
\end{proposition}

\pf
If $A,B\subseteq GF(p)$ then the discussion in \S \ref{sec:1}
implies
\begin{equation}
\label{*}
\wt(r_A r_B)=\sum_{k\in GF(p)} \ {\rm parity}\ |A\cap (k-B)|, \ \ \ \ 
\end{equation}
where $k-B=\lbrace k-b\ |\ b\in B\rbrace$ and parity$(x)=1$ if
$x$ is an odd integer, and $=0$ otherwise. 
Let $S\subseteq GF(p)$, then we have
\[
p-\wt(r_Q r_S)-\wt(r_N r_S)
=\sum_{a\in GF(p)} 
\big{(}
1-\ {\rm parity}\ |Q\cap (a-S)|-\ {\rm parity}\ |N\cap (a-S)|
\big{)}.
\]
Let
\[
T_a(S)=1-\ {\rm parity}\ |Q\cap (a-S)|-\ {\rm parity}\ |N\cap (a-S)|.
\]

{\bf Case 1.} 
If $|S|$ is even and $a\in S$ then $0\in a-S$ so $|Q\cap (a-S)|$ 
odd implies that $|N\cap (a-S)|$ is even, since $0$ is not included 
in $Q\cap (a-S)$ or $N\cap (a-S)$.  Likewise, $|Q\cap (a-S)|$ even 
implies that $|N\cap (a-S)|$ is odd.  Therefore $T_a(S)=0$.

{\bf Case 2.} If $|S|$ is even and $a\notin S$ then parity
$|Q\cap (a-S)|=$parity$|N\cap (a-S)|$.  If $|Q\cap (a-S)|$ is even then 
$T_a(S)=1$ and if $|Q\cap (a-S)|$ is odd then $T_a(S)=-1$.

{\bf Case 3.} $|S|$ is odd.  We claim that $(a-S)^c=a-S^c$.  (Proof: Let 
$s\in S$ and $\bar s\in S^c$.  Then $a-s=a-\bar s\implies s=\bar s$, which 
is obviously a contradiction. Therefore $(a-S)\cap (a-S^c)=\emptyset$, so 
$(a-S)^c\supseteq(a-S^c).$ Replace $S$ by $S^c$ to prove the claim.)  
Also note that 

\[
\big{(} Q\cap (a-S)\big{)} \sqcup \big{(} Q\cap (a-S^c) \big{)}
=GF(p)\cap Q=Q
\]
has $|Q|=\frac{p-1}{2}$ elements ($\sqcup$ denotes disjoint union).  
So

\[
{\rm parity}\ |Q\cap (a-S)|={\rm parity}\ |Q\cap (a-S^c)|
\]
if and only if $|Q|$ is even and 

\[
{\rm parity}\ |Q\cap (a-S)| \neq {\rm parity}\ |Q\cap (a-S^c)|
\]
if and only if and only if $|Q|$ is odd. 

{\bf Conclusion.} 
\[
|S|\ {\rm even\colon}\ T_a(S)=\prod_{x\in a-S}\left( \frac{x}{p}\right)
\]
\[
|S|\ {\rm odd\ and}\ p\equiv 3\pmod 4\colon\ T_a(S)=-T_a(S^c)
\] 
\[
|S|\ {\rm odd\ and}\ p\equiv 1\pmod 4\colon\ T_a(S)=T_a(S^c)
\]
The relation between $\wt(\mathbf{c})$ and the character sum follows from this.
For the remaining part of the equation, use Proposition \ref{prop:hyperelliptic}.
\qed

\begin{remark}
It can be shown, using the coding-theoretic results above, that 
if $p\equiv -1 \pmod 8$ then (for non-empty $S$)
$X_{S}(GF(p))$ contains at least $\sqrt{p}+1$ points. 
This also follows from Weil's estimate, but since the proof is short, it is given 
here.

%\begin{corollary}
%\label{cor:lowerbound}
%If $p\equiv -1 \pmod 8$ and $S$ is non-empty then
%$X_S(GF(p))$ contains at least $\sqrt{p}+2$ points.
%\end{corollary}

%\pf
%This proves the statement of the
%corollary when $|S|$ is even. 
What part (c) of Proposition \ref{prop:main} gives is that
if $p\equiv -1 \pmod 8$ and $|S|$ is odd then
$X_S(GF(p))$ contains at least $\sqrt{p}+2$ points.
If $|S|$ is even then 
perform the substitution 
$x = a+1/\overline{x}$, $y=\overline{y}/\overline{x}^{|S|}$
on the equation $y^2=f_S(x)$. 
This creates a hyperelliptic curve $X$ in $(\overline{x},\overline{y})$ 
for which $|X(GF(p))|=|X_S(GF(p))|$ and $X\cong X_{S'}$, where $|S'|=|S|-1$ is
odd. Now apply part (c) of the above proposition and Remark \ref{remark:sqrt}
to $X_{S'}$. \qed

\end{remark}

\begin{remark}
If $|S|=2$ or $|S|=3$ then more can be said about the character sums 
above.

If $|S|=2$ then $\sum_a \chi(f_S(a))$ can be computed explicitly 
(it is ``usually'' equal to $-1$ - see Proposition 1 in \cite{Wa}).
If $|S|=3$ then $\sum_a \chi(f_S(a))$ can be expressed in terms
of a hypergeometric function $_2F_1$ over $GF(p)$ (see Proposition 2 in 
\cite{Wa}). 

\end{remark}

It has already been observed that the following fact is true.
Since its proof using basic facts about hyperelliptic 
curves is so short, it is included here.

\begin{corollary}
\label{cor:1}
$C_{NQ}$ is an even weight code.
\end{corollary}

\pf 
Since $p$ is odd $1\not= -1$ in $GF(p)$, so every affine
point in $X_S(GF(p))$ occurs as an element of a pair
of solutions of $y^2=f_S(x)$. There are two points at infinity 
(if ramified, it is counted with multiplicity two), so in general 
$|X_S(GF(p))|$ is even. The formulas for the weight of a codeword in 
the above Proposition imply every codeword has even weight.
\qed

As a consequence of this Proposition and Lemma \ref{lemma:mindist},
we have the following result.

\begin{corollary}
\label{cor:5/3bound}
Assume Conjecture \ref{conj:CNQ}.
If $p\equiv 3 \pmod 4$ then
$\max_S\, |X_S(GF(p))|>\frac{5}{3}p -4$.
\end{corollary}

\begin{example}
The following examples were computed with the help of {\tt SAGE}.

If $p=11$ and $S=\{1,2,3,4\}$ then

\[
\begin{array}{c}
(r_S(x)r_N(x),r_S(x)r_Q(x))\\
=(x^{10} + x^9 + x^7 + x^6 + x^5 + x^4 + x^2 + 1,
 x^{10} + x^9 + x^7 + x^6 + x^5 + x^3 + x + 1),
\end{array}
\]
corresponds to the codeword 
$(1, 0, 1, 0, 1, 1, 1, 1, 0, 1, 1, 1, 1, 0, 1, 0, 1, 1, 1, 0, 1, 1)$
of weight $16$. An explicit computation shows that the 
character sum $\sum_{a\in GF(11)} \chi(f_S(a))$ is
$-5$, as expected.

If $p=11$ and $S=\{1,2,3\}$ then

\[
(r_S(x)r_N(x),r_S(x)r_Q(x))
=(x^9 + x^7 + x^5 + x^4 + x^3 + x^2 + x, 
x^{10} + x^8 + x^6 + x^3 + x^2 + x + 1).
\]
corresponds to the codeword 
$(0, 1, 1, 1, 1, 1, 0, 1, 0, 1, 0, 1, 1, 1, 1, 0, 0, 1, 0, 1, 0, 1)$
of weight $14$. An explicit computation shows that the 
character sum $\sum_{a\in GF(11)} \chi(f_{S^c}(a))$ is
$3$, as predicted.

\end{example}

Recall $B(c,p)$ is the statement: 
$|X_S(GF(p))| \leq  c\cdot p$ for all $S\subset GF(p)$.

\begin{theorem}
(Bazzi-Mitter) 
Fix $c\in (0,2)$. If $B(c,p)$ holds for infinitely many $p$
with $p\equiv 1\pmod 4$
then there exists an infinite family of binary codes with asymptotic
rate $R = 1/2$ and relative distance $\delta \geq 1-\frac{c}{2}$.
\end{theorem}

This is an easy consequence of the above Proposition and is
essentially in \cite{BM} (though they assume $p\equiv 3\pmod 8$).
%%
%% BM has a Prop which allows them, at a small cost, to assume
%% |S| is even. I don't understand their proof. If you believe this
%% Prop of theirs then youcan omit ``with $p\equiv 1\pmod 4$''

\begin{theorem}
\label{thrm:main}
If $B(1.77,p)$ is true for infinitely many primes
$p$ with $p\equiv 1\pmod 4$ then Goppa's conjecture is false.
\end{theorem}

\pf
Recall Goppa's conjecture is that the binary asymptotic 
Gilbert-Varshamov bound is best possible for any family of binary codes.  
The asymptotic GV bound states that the rate $R$ 
is greater than or equal to $1-H_2(\delta)$, where 

\[
H_q(\delta)=\delta\cdot \log_q(q-1)-
\delta\log_q(\delta)-(1-\delta)\log_q(1-\delta)
\]  
is the entropy function (for a $q$-ary channel).
Therefore, according to Goppa's conjecture, 
if $R=\frac{1}{2}$ (and $q=2$) then the best possible $\delta$ is 
$\delta_0=.11$.  Assume $p\equiv 1\pmod 4$.
Goppa's conjecture implies that the minimum distance of our 
QQR code with rate $R=\frac{1}{2}$ satisfies $d<\delta_0\cdot 2p=.22p$,
for sufficiently large $p$.
Recall that the weight of a codeword in this QQR code 
is given by Proposition \ref{prop:main}.
$B(1.77,p)$ (with $p\equiv 1\pmod 4$) implies (for all $S\subset GF(p)$)
$\wt((r_Sr_N,r_Sr_Q))\geq 2p-|X_S(GF(p))|\geq 0.23p$. In other words, for
$p\equiv 1\pmod 4$, all nonzero codewords have weight at least
$0.23p$.
This contradicts the estimate above.\qed

Using the same argument and the 
first McEliese-Rumsey-Rodemich-Welsh (MRRW) bound (\cite{HP}, Theorem 2.10.6),
we prove the following unconditional result.

\begin{theorem}
For all sufficiently large primes $p$ for which $p\equiv 1\pmod 4$,
the statement $B(1.62,p)$ is false.
\end{theorem}

\pf
If a prime $p$ satisfies $B(1.62,p)$ then we shall call it ``admissible.''
We show that the statement ``$B(1.62,p)$ holds for all sufficiently large
primes $p$ for which $p\equiv 1\pmod 4$''
contradicts the first asymptotic MRRW bound.
Indeed, this MRRW bound states that the rate $R$ 
is less than or equal to 

\[
h(\delta)=H_2(\frac{1}{2} - \sqrt{\delta(1-\delta)}).
\]
This, and the fact that $R=\frac{1}{2}$ for our QQR codes 
(with $p\equiv 1\pmod 4$), imply
$\delta \leq \delta_0 = h^{-1}(1/2)\cong 0.187$. 
Therefore, for all large $p$ (admissible or not),
$d\leq \delta_0\cdot 2p$. On the other hand, if $p$ is admissible and
$|X_S(GF(p))| \leq c\cdot p$ (where $c=1.62$) then 
by the above argument, $d\geq 2\cdot(p-\frac{c}{2}p)$.
Together, we obtain $1-\frac{c}{2}\leq\delta_0$, so
$c\geq 2\cdot (1-h^{-1}(1/2))\cong 1.626$. This is a contradiction.
\qed

\begin{corollary}
\label{cor:main}
Assume Conjecture \ref{conj:CNQ}.
There is a constant $p_0$ (ineffectively computable) having the following
property: if $p>p_0$ then there is a 
subset $S\subset GF(p)$ for which the bound $|X_S(GF(p))| > 1.62p$
holds.
\end{corollary}

This is of course the same as the above theorem, except that
we have used Corollary \ref{cor:5/3bound} (which 
unfortunately depends on Conjecture \ref{conj:CNQ})
to remove the hypothesis $p\equiv 1\pmod 4$.

\section{Weight distributions}
\label{sec:duursma}

In \cite{D1} Iwan Duursma associates to a linear code $C$ 
over $GF(q)$ a {\it zeta function} $Z=Z_C$ of the form

\[
Z(T)=\frac{P(T)}{(1-T)(1-qT)},
\] 
where $P(T)$ is a polynomial of degree $n+2-d-d^\perp$
which only depends on $C$ through its weight enumerator 
polynomial
(here $d$ is the minimum distance of $C$ and $d^\perp$ is
the minimum distance of its dual code $C^\perp$; we assume
$d\geq 2$ and $d^\perp\geq 2$).
If $\gamma=\gamma(C) = n+k+1-d$ and 
$z_C(T)=Z_C(T)T^{1-\gamma}$ then 
the functional equation in \cite{D1} can be written in the form
$z_{C^\perp}(T)=z_C(1/qT)$. If we let
$\zeta_C(s) = Z_C(q^{-s})$ and $\xi_C(s) = z_C(q^{-s})$
then $\zeta_C$ and $\xi_C$ have the same zeros but $\xi_C$ is
``more symmetric'' since the functional equation expressed in terms of it
becomes

\[
\xi_{C^\perp}(s)=\xi_{C}(1-s).
\]
Abusing terminology, we call both $Z_C$ and $\zeta_C$ a
{\it Duursma zeta function}.
In fact, if $\rho_i$ denotes the $i$-th zero of the zeta function
$Z(T)$ of an actual code then equations (5)-(6) of \cite{D2} implies 
(for the even weight binary codes we are considering here) the relation 

\[
d= 2-\sum_{i} \rho_i^{-1}. 
\]
Therefore, further 
knowledge of the zeros of $Z(T)$ could be very useful.

If $C$ is self-dual (or actually only formally self-dual)
then the zeros of the $\zeta$-function occur in 
pairs about the ``critical line'' $Re(s)=\frac{1}{2}$.
Following Duursma, we say (for formally self-dual codes
$C$) the zeta function $\zeta_C$ satisfies the
{\it Riemann hypothesis} if all its zeros occur on the ``critical line''.

\begin{example}
The following computations were done with the help of {\tt SAGE}.
If $p=7$ then 
the $[14,7,4]$ (self-dual) code $C_{NQ}$ has ``zeta polynomial''

\[
P(T)= \frac{2}{143} + \frac{4}{143}T + \frac{19}{429}T^2 + 
\frac{28}{429}T^3 + \frac{40}{429}T^4 + \frac{56}{429}T^5 + 
\frac{76}{429}T^6 + \frac{32}{143}T^7 + \frac{32}{143}T^8.
\]
It can be checked that all the roots $\rho$ of $Z_C$
have $|\rho|=1/\sqrt{2}$, thus verifying the Riemann hypothesis in this case.
\end{example}

It would be interesting to know if the Duursma zeta function 
$Z(T)$ of $C_{NQ}$, for $p\equiv 3 \pmod 4$, always satisfies the 
Riemann hypothesis.

A self-dual code is called
{\it extremal} if its minimum distance satisfies the Sloane-Mallows bound
\cite{D3} and {\it optimal} if its minimum distance is maximal among all such
linear codes of that length and dimension (see also Chinen 
\cite{Ch}, \cite{Ch2}). 
As noted above, the Duursma zeta function only depends on the weight
enumerator. It has been conjectured that, for all extremal self-dual codes $C$, 
the $\zeta$-function satisfies the Riemann hypothesis.
The example below shows that
``extremal self-dual'' cannot be replaced by ``optimal formally self-dual''.

Based on computer computations using \sage, the following 
statement appears to be true, though we have no proof.

\begin{conjecture}
If $p\equiv 1 \pmod 4$ then the code
$C'$ spanned by $C_{NQ}$ and the all ones codeword
(i.e., the smallest code containing $C_{NQ}$ and all its
complementary codewords) is a formally self-dual code
of dimension $p$. Moreover, we if $A=[A_0,A_1,...,A_n]$ denotes
the weight distribution vector of $C_{NQ}$ then
the weight distribution vector of 
$C'$ is $A+A^*$, where $A^*=[A_n,...,A_1,A_0]$. 
\end{conjecture}

Using {\tt SAGE}, it can be shown that the Riemann hypothesis is not valid 
for these ``extended QQR codes'' in general, as the 
following example illustrates.

\begin{example}
If $p=13$ then $C'$ is a $[26,13,6]$ code with
weight distribution

\[
[1, 0, 0, 0, 0, 0, 39,
 0, 455, 0, 1196, 0, 2405,
 0, 2405, 0, 1196, 0, 455,
 0, 39, 0, 0, 0, 0, 0, 1].
\]
This is (by coding theory tables, as included in {\tt SAGE}
\cite{S}) an optimal, formally self-dual code.
This code $C'$ has zeta polynomial

%{\small{
\[
\begin{array}{ll}
P(T)&= \frac{3}{17710} + \frac{6}{8855}T + \frac{611}{336490}T^2 + 
\frac{9}{2185}T^3 + \frac{3441}{408595}T^4 + 
\frac{6448}{408595}T^5 + \frac{44499}{1634380}T^6  \\
&+
\frac{22539}{520030}T^7 + \frac{66303}{1040060}T^8 + 
\frac{22539}{260015}T^9 + \frac{44499}{408595}T^{10} + 
\frac{51584}{408595}T^{11}  \\
& +\frac{55056}{408595}T^{12} + 
\frac{288}{2185}T^{13} + \frac{19552}{168245}T^{14} + 
\frac{768}{8855}T^{15} + \frac{384}{8855}T^{16}.
\end{array}
\]
%}}
Using {\tt SAGE}, it can be checked that 
only 8 of the 12 zeros of this function
have absolute value $1/\sqrt{2}$.
\end{example}

\section{Long Quadratic Residue Codes}
\label{sec:lqrc}
\label{sec:LQR}

We now introduce a new code, constructed similarly to the
QQR codes discussed above:

\[
C=\lbrace ({r}_N r_S, {r}_Q r_S, {r}_N r_S^*, {r}_Q r_S^*)\ |\
S\subseteq GF(p)\rbrace .
\]
We call this a {\it long quadratic residue code} or
{\it LQR code} for short, and identify it with a subset of
$\fff^{4p}$.
Observe that this code is non-linear.

For any $S\subseteq GF(p)$, let

\[
\mathbf{c}_S=({r}_N r_S, {r}_Q r_S, {r}_N r_S^*, {r}_Q r_S^*)
\]
and let

\[
v_S=({r}_N r_S, {r}_Q r_S, {r}_N r_{S}, {r}_Q r_{S}).
\]
If $S_1\Delta S_2$ denotes the symmetric difference between
$S_1$ and $S_2$ then it is easy to check that

\begin{equation}
\label{eqn:add}
\mathbf{c}_{S_1}+\mathbf{c}_{S_2}=v_{S_1\Delta S_2}.
\end{equation}

We now compute the size of $C$ using Lemma \ref{lemma:prod}. 
We prove the \textit{claim}:
if $p\equiv 3\pmod4$ then the map that sends $S$ to the codeword
$\mathbf{c}_S$ is injective. 
This implies $|C|=2^p$.  Suppose not, then there are two subsets
$S_1,\ S_2\subseteq GF(p)$ that are mapped to the same codeword.
Subtracting $\mathbf{c}_{S_1}-\mathbf{c}_{S_2}=\mathbf{c}_{S_1}+\mathbf{c}_{S_2}=v_{S_1\Delta S_2}$, and
the subset $T=S_1\Delta S_2$ satisfies
$r_Q r_T=r_N r_T=r_Q r_{T^c}=r_N r_{T^c}=0$. If $|T|$ is even then
$0=(r_Q + r_N)r_T=(r_{GF(p)}-1)r_T=r_T$.  This forces $T$ to be
the empty set, so $S_1=S_2$.  Now if $|T|$ is odd then similar
reasoning implies that $T^c$ is the empty set.  Therefore,
$S_1=\emptyset$ and $S_2=GF(p)$ or vice versa.  This proves the claim.

In case $p\equiv 1\pmod4$, we \textit{claim}: $|C|=2^{p-1}$.
Again, suppose there are two subsets
$S_1,\ S_2\subseteq GF(p)$ that are mapped to the same codeword.
Then the subset $T=S_1\Delta S_2$ satisfies
$r_Q r_T=r_N r_T=r_Q r_{T^c}=r_N r_{T^c}=0$. This implies
either $T=\emptyset$ or $T=GF(p)$. Therefore, either
$S_1=S_2$ or $S_1=S_2^c$.

Combining this discussion with Proposition \ref{prop:hyperelliptic},
we have proven the following result.

\begin{theorem}
\label{thrm:1}
The code $C$ has length $n=4p$ and has size $M=2^{p-1}$
if $p\equiv 1\pmod 4$,
and size $M=2^p$ if $p\equiv 3\pmod 4$.
If $p\equiv 3\pmod 4$ then the minimum non-zero weight is $2p$
and the minimum distance is at least

\[
d_p=4p-2\max_{S\subset GF(p)} |X_S(GF(p))|.
\]
If $p\equiv 1\pmod 4$ then $C$ is a binary
$[4p,p-1,d_p]$-code.
\end{theorem}

\begin{remark}
If $p\equiv 3\pmod 4$, there is no simple reason I can think of 
why the minimum distance should actually be less than the minimum non-zero 
weight.
\end{remark}

\begin{lemma}
If $p\equiv 1\pmod 4$ then 

\begin{itemize}
\item
$v_S=\mathbf{c}_S$, 

\item
$\mathbf{c}_{S_1}+\mathbf{c}_{S_2}=\mathbf{c}_{S_1\Delta S_2}$,
\item
the code $C$ is isomorphic to the QQR code $C_{NQ}$.

\end{itemize}
\end{lemma}

In particular, $C$ is linear and of dimension $p-1$.

\pf
It follows from the the proof of Theorem \ref{thrm:1}
that if $p\equiv 1\pmod 4$ then
$r_Nr_{S_1}=r_Nr_{S_2}$ and $r_Qr_{S_1}=r_Qr_{S_2}$
if and only if $S_2=S_1^c$. The lemma follows rather easily
as a consquence of this and (\ref{eqn:add}).
\qed

Assume $p\equiv 3\pmod 4$. Let

\[
V = \{v_S\ |\ S\subset GF(p)\}
\]
and let 

\[
\overline{C}=C\cup V.
\]

\begin{lemma}
\label{lemma:p3mod4}
The code $\overline{C}$ is 

\begin{itemize}
\item
the smallest linear subcode of $\fff^{4p}$ containing $C$,
\item
dimension $p+1$,
\item
minimum distance $\min(d_p,2p)$.
\end{itemize}
\end{lemma}

By abuse of terminology, we call $\overline{C}$ an {\it LQR code}.

\pf
The first part follows from (\ref{eqn:add}). The second part follows 
from a counting argument (as in the proof of Theorem \ref{thrm:1}).
The third part is a corollary of Theorem \ref{thrm:1}.
\qed

Recall that

\[
\wt({r}_N r_S,{r}_Q r_S)
=
\left\{
\begin{array}{ll}
p-\sum_{a\in GF(p)} \left(\frac{f_S (a)}{p}\right),
&|S|\ {\rm even\ (any\ }p {\rm )}, \\
p-\sum_{a\in GF(p)} \left(\frac{f_{S^c} (a)}{p}\right), 
& \ |S|\ {\rm odd\ and\ }p\ \equiv 1\pmod4,\\
p+\sum_{a\in GF(p)} \left(\frac{f_{S^c} (a)}{p}\right),
&|S|\ {\rm odd\ and}\ p\equiv 3\pmod 4,
\end{array}
\right.
\]
by Proposition \ref{prop:main}. 

\begin{lemma}
\label{lemma:wtLQR}
For each $p$, the codeword
$\mathbf{c}_S=\left( {r}_N r_S, {r}_Q r_S, {r}_N r_S^*, {r}_Q r_S^*\right)$
of $C$ has weight

\[
\wt(\mathbf{c}_S)
=
\left\{
\begin{array}{ll}
2p-2\sum_{a\in GF(p)} \left(\frac{f_S (a)}{p}\right),& p\equiv 1\pmod 4,\\
2p, & p\equiv 3\pmod 4.
\end{array}
\right.
\]
\end{lemma}

In other words, if $p\equiv 3\pmod 4$ then $C$ is a constant weight code.

\pf
Indeed, Proposition \ref{prop:main} implies
if $p\equiv 1\pmod 4$ then 

\begin{equation}
\label{eqn:symm}
\begin{array}{ll}
\wt\left( {r}_N r_S, {r}_Q r_S, {r}_N r_S^*, {r}_Q r_S^*\right)
&=\ \  \wt\left( {r}_N r_S, {r}_Q r_S\right)+
\wt\left({r}_N r_S^*, {r}_Q r_S^*\right)\\
&=\ \  2\cdot \wt\left( {r}_N r_S, {r}_Q r_S\right)\\
&=\ \  2p-2\sum_{a\in GF(p)} \left(\frac{f_S (a)}{p}\right),
\end{array}
\end{equation}
if $p\equiv 3\pmod 4$ and $|S|$ is even then

\begin{equation}
\label{eqn:symm3}
\begin{array}{ll}
\wt\left( {r}_N r_S, {r}_Q r_S, {r}_N r_S^*, {r}_Q r_S^*\right)
&=\ \  \wt\left( {r}_N r_S, {r}_Q r_S\right)+
\wt\left({r}_N r_S^*, {r}_Q r_S^*\right)\\
&=\ \  p-\sum_{a\in GF(p)} \left(\frac{f_{S} (a)}{p}\right)+
p+\sum_{a\in GF(p)} \left(\frac{f_{S} (a)}{p}\right)\\
&=\ \  2p,
\end{array}
\end{equation}
and if $p\equiv 3\pmod 4$ and $|S|$ is odd then

\begin{equation}
\label{eqn:symm4}
\begin{array}{ll}
\wt\left( {r}_N r_S, {r}_Q r_S, {r}_N r_S^*, {r}_Q r_S^*\right)
&=\ \  \wt\left( {r}_N r_S, {r}_Q r_S\right)+
\wt\left({r}_N r_S^*, {r}_Q r_S^*\right)\\
&=\ \  p+\sum_{a\in GF(p)} \left(\frac{f_{S} (a)}{p}\right)+
p-\sum_{a\in GF(p)} \left(\frac{f_{S} (a)}{p}\right)\\
&=\ \  2p.
\end{array}
\end{equation}
\qed

\begin{example}
The following examples were computed with the help of {\tt SAGE}.
When $p=11$ and $S=\{1,2,3,4\}$, $\mathbf{c}_S$
%
%{\tiny{
%\[
%\begin{array}{l}
%\left( {r}_N r_S, {r}_Q r_S, {r}_N r_S^*, {r}_Q r_S^*\right)=
%\\
%(x^9 + x^7 + x^5 + x^4 + x^3 + x^2 + x,
% x^{10} + x^8 + x^6 + x^3 + x^2 + x + 1,
% x^{10} + x^8 + x^6 + 1,
% x^9 + x^7 + x^5 + x^4)
%\end{array}
%\]
%}}
corresponds to the codeword

{\small{
\[
(0, 1, 1, 1, 1, 1, 0, 1, 0, 1, 0,
1, 1, 1, 1, 0, 0, 1, 0, 1, 0, 1,
1, 0, 0, 0, 0, 0, 1, 0, 1, 0, 1,
0, 0, 0, 0, 1, 1, 0, 1, 0, 1, 0)
\]
}}
of weight $22$. When $p=11$ and $S=\{1,2,3\}$, $\mathbf{c}_S$
%
%{\tiny{
%\[
%\begin{array}{c}
%\left( {r}_N r_S, {r}_Q r_S, {r}_N r_S^*, {r}_Q r_S^*\right)
%\\
%=(x^{10} + x^9 + x^7 + x^6 + x^5 + x^4 + x^2 + 1,
% x^{10} + x^9 + x^7 + x^6 + x^5 + x^3 + x + 1,
% x^8 + x^3 + x,
% x^8 + x^4 + x^2)
%\end{array}
%\]
%}}
corresponds to the codeword

{\small{
\[
(1, 0, 1, 0, 1, 1, 1, 1, 0, 1, 1,
 1, 1, 0, 1, 0, 1, 1, 1, 0, 1, 1,
 0, 1, 0, 1, 0, 0, 0, 0, 1, 0, 0,
 0, 0, 1, 0, 1, 0, 0, 0, 1, 0, 0)
\]
}}
of weight $22$.

\end{example}

It turns out Lemma \ref{lemma:p3mod4} allows us to 
improve the statement of Theorem \ref{thrm:main} in \S \ref{sec:BM}. 
The next subsection is devoted to this goal.

\subsection{Goppa's conjecture revisited}

We shall now remove the condition $p\equiv 1 \pmod 4$ in one of the
results in \S \ref{sec:BM}, 
at a cost of weakening the constant involved.

Assuming $B(c,p)$ holds, we have that the minimum distance of
$\overline{C}$ is $\geq \min(d_p,2p) \geq 4p(1-\frac{c}{2})$
and the information rate is $R=\frac{1}{4}+\frac{1}{4p}$.
When $R=1/4$, Goppa's conjecture gives $\delta = 0.214...$\, . 
So Goppa's conjecture will be 
false if $1-\frac{c}{2}=0.215$, or $c=1.57$.
We have the following improvement of Theorem \ref{thrm:main}.

\begin{theorem}
\label{thrm:impr}
If the $B(1.57,p)$ is true for infinitely many primes $p$ then 
Goppa's conjecture is false.
\end{theorem}

A similar argument (using $h(x)$ and the MRRW bound in place of
$1-H_2(x)$ and the hypothetical Goppa bound) gives

\begin{theorem}
\label{thrm:lqr}
$B(1.39,p)$ cannot be true for infinitely many primes $p$.
In other words, for all ``sufficiently large'' $p$,
we must have $X_S(GF(p))>1.39p$ for some $S\subset GF(p)$.
\end{theorem}

\section{Some results of Voloch}
\label{sec:voloch}

\begin{lemma} (Voloch)
If $p\equiv 1,3 \pmod 8$ then $|X_Q(GF(p))| = 1.5p+a$, where $Q$ is 
the set of quadratic residues and $a$ is a small
constant, $-\frac{1}{2}\leq a \leq \frac{5}{2}$.
\end{lemma}

A similar bound holds if $X_Q$ is replaced by $X_N$ and 
$p\equiv 1,3 \pmod 8$ is replaced by $p\equiv 7 \pmod 8$
(in which case $2$ is a quadratic residue).

\vskip .1in

\pf
By Proposition \ref{prop:hyperelliptic}, we know
that if $p\equiv 3 \pmod 8$
(so $|Q|$ is odd):

\[
\sum_{a\in GF(p)} \chi(f_Q(a))
=-p-1+|X_Q(GF(p))|.
\]
Similarly,  if $p\equiv 1 \pmod 8$
(so $|Q|$ is even):

\[
\sum_{a\in GF(p)} \chi(f_Q(a))
=-p-2+|X_Q(GF(p))|.
\]

Since $b^{\frac{p-1}{2}} \equiv \chi(b) \pmod p$, we have
 
\[ 
x^{\frac{p-1}{2}}-1 = \prod_{a\in Q} (x-a)=f_Q(x),
\ \ \ \ 
x^{\frac{p-1}{2}}+1 = \prod_{a\in N} (x-a).
\]
In particular, for all $n\in N$,

\[
f_Q(n) = \prod_{a\in Q} (n-a) = n^{\frac{p-1}{2}}-1 \equiv -2 \pmod p.
\]
Since $p\equiv 1,3 \pmod 8$, we have $\chi(-2)=1$,
so $\chi(f_Q(n))=1$ for all $n\in N$. It follows that
$|X_Q(GF(p))|= \frac{3}{2}p+\chi(f_Q(0))+\frac{1}{2}$
(if $p\equiv 3 \pmod 8$) or 
$|X_Q(GF(p))|= \frac{3}{2}p+\chi(f_Q(0))+\frac{3}{2}$
(if $p\equiv 1 \pmod 8$). 
\qed

\vskip .1in

Here is an extension of the idea in the above proof. Fix 
an integer $\ell>2$. Assuming $\ell$ divides $p-1$, there are distinct
$\ell$-th roots $r_1=1$, $r_2$, ..., $r_\ell$ in $GF(p)$ for which
$x^{p-1}-1 = \prod_{i=1}^\ell (x^{\frac{p-1}{\ell}}-r_i)$.
Also, $x^{\frac{p-1}{\ell}}-1=\prod_{a\in P_\ell}(x-a)=f_{P_\ell}(x)$,
where $P_\ell$ denotes the set of non-zero $\ell$-th powers in 
$GF(p)$. 

{\it Claim}: It is possible to find an infinite sequence
of primes $p$ satisfying $p\equiv 1 \pmod \ell$ and 
$\chi(r_i-1)=1$, for all $2\leq i\leq \ell$ (where $\chi$ denotes the
Legendre character mod $p$).
If the claim is true then we will have a lower
bound for $|X_{P_\ell}(GF(p))|$ on the order of $(2-\frac{1}{\ell})p$,
along the lines above, by Proposition \ref{prop:hyperelliptic}.

Proof of claim: It is a well-known fact in algebraic number theory 
that $p\equiv 1 \pmod \ell$ implies that the prime $p$
splits completely in the cyclotomic field $\qqq_\ell$ generated by
the $\ell$-th roots of unity in $\ccc$, denoted $\tilde{r}_1=1$, $\tilde{r}_2$,
..., $\tilde{r}_\ell$. The condition $\chi(r_i - 1) = 1$ means that $p$ 
splits in the extension of $\qqq_\ell$ obtained by adjoining 
$\sqrt{\tilde{r}_i-1}$ (here $i=2,...,\ell$). By Chebotarev's
density theorem there exist infinitely many such $p$, as claimed. 
\qed

In fact, there are effective versions which give explicit 
information on computing such $p$ \cite{LO}, \cite{Se}.
This, together with the previous lemma, proves the following result.

\begin{theorem} (Voloch)
If $\ell\geq 2$ is any fixed integer then for infinitely many
primes $p$ there exists a subset $S\subset GF(p)$ for which
$|X_S(GF(p))| = (2-\frac{1}{\ell})p+a$, where $a$ is a small
constant, $-\frac{1}{2}\leq a \leq \frac{5}{2}$.
\end{theorem}
 
In fact, the primes occurs with a positive (Dirichlet) density
and the set $S$ can be effectively constructed.

\vskip .4in
{\it Acknowledgements}: I thank Prof. Amin Shokrollahi of the
Ecole Polytechnique F\'ed\'erale de Lausanne for helpful advice
and Prof. Felipe Voloch  of the University of Texas 
for allowing his construction to be included above.
Parts of this note (such as Proposition \ref{prop:main})
can be found in the honors thesis \cite{C} of my former 
student Greg Coy, who was a pleasure to work with.
Last but not least, I thank the anonymous referees
who substantially improved the paper with their 
insightful suggestions.

\end{document}